\documentclass[12pt,draftcls,onecolumn]{IEEEtran}
\usepackage{cite}
\usepackage{graphicx}
\usepackage{subfig}
\usepackage{pgf} 
\usepackage{bbold} 
\usepackage{amsmath}

\begin{document}
\title{Laplace Transform of Product of Generalized Marcum Q, Bessel I, and Power Functions with Applications\footnote{This work has been submitted to IEEE for possible publication. Copyright may be transferred without notice, after which this version may no longer be accessible.
} }
\author{
Natalia Y.~Ermolova,
 and Olav Tirkkonen,~\IEEEmembership{Member,~IEEE}\\
\small{Department of Communications and Networking, Aalto University, P.O. Box 13000, FI-00076 Aalto,   
  Finland,\\
 e-mail: natalia.ermolova@aalto.fi; olav.tirkkonen@aalto.fi}
}
\maketitle
\begin{abstract}
 The evaluation of integral transforms of special functions is required in different research and practical areas.  Analyzing  the $\kappa$-$\mu$ fading distribution also khown as the generalized Rician distribution, we find out that the assessment of a few different performance metrics, such as the  probability of energy detection of unknown signals and outage probability under co-channel interference, involves the evaluation of an infinite-range integral, which has the form of Laplace transform of product of Marcum Q, Bessel I, and power functions. We evaluate this integral in a closed-form and present numerical estimates.
\end{abstract}
\section{Introduction} The evaluation of integral transforms of special functions is required in different areas of engineering. This problem arises when solving  various analysis and design tasks such as assessment of average performance metrics. In some practical applications,  integral transforms of 
 the Marcum Q and Bessel I functions occur. In communication engineering, these functions appear in statistical models of fading radio channels, see, for example,  [1]-[4].  The Marcum Q function arises also in different tasks of signal detection in additive white Gaussian noise (AWGN) [4]-[8].

Both functions are available via standard software packages such as \emph{Mathematica}, but nowadays analytical results on integral transforms of these special functions are not widely presented in the literature, and only some special cases have been reported. A large number of integrals involving the Marcum Q function of first order is presented in\cite{Nut}.  

Due to the modern software, even infinite-range integrals often can be evaluated numerically. But it is well known that this process may be accompanied by numerical problems. Therefore, it is convenient to have closed-form expressions, which additionally can be useful at the design stage since they allow analyzing of  effects of various system parameters on performance metrics. Furthermore, closed-form expressions are useful for solving various optimization tasks. 

The product of Bessel I and power functions  
 occurs in the probability density function (PDF) of the $\kappa$-$\mu$ generalized fading distribution recently introduced by M.-D. Yacoub  for modeling propagation effects  in a non-homogeneous line-of-sight environment \cite{Yac}. The $\kappa$-$\mu$ distribution is also known as the generalized Rician distribution [9]-[10]. Analyzing this fading distribution, we reveal that the same integral expressed as the Laplace transform of product of  Marcum Q, Bessel I, and power functions occurs in a few 
practical tasks related to the performance assessment of communication systems over $\kappa$-$\mu$ fading. Such are evaluation of probability of energy detection of unknown signals as well as assessment  of outage probability 
in $\kappa$-$\mu /\kappa$-$\mu$ interference-limited scenarios. As far as we aware, both  tasks have not been yet solved analytically.

The detection probability over $\kappa$-$\mu$ fading was evaluated so far by applying series expansions of Marcum Q and Bessel I functions 
[11]-[12]. These procedures resulted in formulas given in terms of infinite series \cite{Annam} that can be expressed in terms of multivariate hypergeometric functions \cite{Sof}. The hypergeometric functions required for the evaluation are not nowadays available via the standard software.

Considering the $\kappa$-$\mu$ fading channels with $\kappa$-$\mu$--faded co-channel interference, we note that no analytical results on the outage probability have been reported so far, while other generalized fading scenarios, such as $\eta$-$\mu$/$\eta$-$\mu$,  $\eta$-$\mu$/$\kappa$-$\mu$, and $\kappa$-$\mu$/$\eta$-$\mu$ cases, have been successfully analyzed [13]-[15].

In this paper, we evaluate the Laplace transform ($\mathcal{L}$) of the product of  the generalized Marcum Q,  Bessel I, and power functions. More precisely, we solve the integral of the form
\begin{IEEEeqnarray*}{c}
\label{int1}
\mathrm{In}(\alpha,\beta,c,p,\mu _1,\mu _2)=\mathcal{L}\left\{Q_{\mu _1}(\alpha \sqrt{t},\beta)t^{\frac{\mu _2-1}{2}}I_{\mu_2 -1}(c\sqrt{t}),\left\{t,p\right\}\right\}\\
=\int_0^{\infty}\mathrm{exp}(-pt)Q_{\mu _1}(\alpha \sqrt{t},\beta)t^{\frac{\mu _2-1}{2}}I_{\mu_2 -1}(c\sqrt{t})dt\\
=2\int_0^{\infty}\mathrm{exp}(-pt^2)t^{\mu _2}Q_{\mu _1}(\alpha t,\beta)I_{\mu_2 -1}
(ct)dt
\IEEEyesnumber
\end{IEEEeqnarray*}
where $I_\nu (.)$ is the modified Bessel function of the first kind and order $\nu$ \cite{Bat}, $Q_M (\alpha ,\beta)$ is the generalized Marcum Q function [4]-[8], [16]-[17] $\mu _1$, $\mu _2$, $c$$\in (-\infty, \infty)$, $p$ $\in \left(0, \infty\right)$, and $\alpha$, $\beta$ are either real or purely imaginary \cite{Bla}. 

We use the derived results  be used for the evaluation of  outage probability over interference-limited $\kappa$-$\mu$ fading radio channels with $\kappa$-$\mu$--faded co-channel interference. We also apply the obtained formulas to assess the detection probability of unknown signals in $\kappa$-$\mu$ fading.

\section{Preliminaries}
In this section, we introduce concepts and special functions used in this work.
\subsection{The $\kappa$-$\mu$ fading distribution}
 The $\kappa$-$\mu$ fading signal is considered as the composition of clusters of multipath waves with the uncorrelated Gaussian in-phase (I) and quadrature (Q) components within each cluster, and
 the probability density function (PDF) $f_{\gamma _{\kappa-\mu}}$  of the $\kappa$-$\mu$ power variable $\gamma _{\kappa -\mu}$ is given in  \cite{Yac} as
\begin{IEEEeqnarray*}{c}
f_{\gamma_{\kappa-\mu}}(x)=\frac{\mu {(1+\kappa)}^{\frac{\mu+1}{2}}x^{\frac{\mu-1}{2}}} 
{\kappa^{\frac{\mu-1}{2}}\mathrm{exp}(\mu \kappa) 
 \Omega ^{\frac{\mu+1}{2}}} \mathrm{exp}\left(-\frac{ \mu (1+\kappa) x }{\Omega}\right)
 I_{\mu-1}\left(2 \mu 
\sqrt{\frac{\kappa (1+\kappa)x}{ \Omega  }}\right)
\IEEEyesnumber
\end{IEEEeqnarray*}
where $\kappa>0$ is the ratio of the total power of the dominant components to that of the scattered waves, $\Omega  =E\{\gamma_{\kappa -\mu}\}$, and 
$\mu=\frac{\Omega _\kappa ^2}{2 \mathrm{var}\{\gamma _{\kappa -\mu}\}} \cdot\frac{
1+2 \kappa}{(1+\kappa)^2}$ ( $E\{\}$ and $\mathrm{var}\{\}$ denote the expectation and variance, respectively)  characterizes the number of multipath clusters. Three parameters of the distribution provide a better fit to experimental data than do the commonly used fading models. As it is pointed out in \cite{Yac}, a fitting procedure may result in non-integer values of $\mu$ that can be caused by a few factors such as the non-Gaussian distribution of the I and Q components, non-zero correlation between the I and Q components, or non-zero correlation between the multipath clusters.

 The  $\kappa$-$\mu$ distribution is a generalized fading model that includes many well-known fading distributions as particular cases. For example, the Rice distribution ($\mu =1$) and Nakagami-$m$ distribution ($\kappa \to 0$) are particular cases of model (2). For integer values of $\mu$, (2) reduces to the non-central chi square distribution.

The cumulative distribution function (CDF) corresponding to (2) is given by \cite[eq. (3)] {Yac} as
\begin{equation}
F_{\gamma_{\kappa-\mu}}(z)=1-Q_{\mu }\left[\sqrt {2\kappa \mu }, \sqrt{\frac{2(1+\kappa )\mu z}{\Omega}}\right].
\end{equation}
\subsection{Modified Bessel function of the first kind (Bessel I function)}
The Bessel I function of the arbitrary order $\nu$ is one of the solution of modified Bessel differential equation  \cite{Bat}, and it is often defined via it series expansion as \cite{Bat}
\begin{equation}
I_{\nu}(x)=\sum_{k=0}^{\infty}\frac{\left(\frac {x}{2}\right)^{2k+\nu}}{\Gamma (\nu +k+1)k!}
\end{equation}
where $\Gamma (.)$ is the gamma function. For integer values of $ \nu =n$, $I_{-n}=I_n$.

\subsection{Confluent hypergeometric functions of two variables $\Phi _3 (b;g;w,z)$ and $\Psi _2(a;d,d';w,z)$}
The confluent hypergeometric functions of two variables  $\Phi _3(b;g;w,z)$ and $\Psi _2(a;d,d';w,z)$ are defined via absolutely convergent hypergeometric series as \cite [vol. 3, eq. (7.2.4.7)]{Prud} and \cite [vol. 3, eq. (7.2.4.9)] {Prud}
\begin{equation}
\Phi _3(b;g;w,z)=\sum_{k,l=0}^{\infty}\frac{(b)_k}{(g)_{k+l}k!l!}w^kz^l,
\end{equation}
and
\begin{equation}
\Psi _2(a;d,d';w,z)=\sum_{k,l=0}^{\infty}\frac{(a)_{k+l}}{(d)_{k}(d')_{l}k!l!}w^kz^l
\end{equation}
where $(a)_k$ means the Pochhammer index \cite [vol. 3, Section II.2] {Prud}, and $-g, -d,-d' \notin \mathbb{Z}^*$ (with $\mathbb{Z}^*$ denoting the set of positive integers ($\mathbb{Z}^+$) and 0, that is $\mathbb{Z}^*=\mathbb{Z}^+\cup \{0\}$).

The Laplace transform of the product $t^{g-1}\Phi _3(b,g,\varsigma t,wt)$ is given in \cite[vol. 4, eq. (3.43.8)]{Prud} as 
\begin{equation}
\mathcal{L}\left\{t^{g-1}\Phi _3(b,g,\varsigma t,wt),\left\{t,p\right\}\right\}=\frac{\Gamma (g)}{p^{g-b}}(p-\varsigma)^{-b}\mathrm{exp}\left(w/p\right)
\end{equation}
where the real parts of $g$, $p$, and $(p-\varsigma )$ are positive.

A series expansion of $\Phi _3(b;g;w,z)$ is given in \cite[eq. (29)]{Par2} as
\begin{IEEEeqnarray*}{c}
\Phi _3(b;g;z,w)=\Gamma (g)w^{\frac{1-g}{2}}\sum_{j=0}^{\infty}\frac{(b)_j}{j!}\left(\frac{z}{\sqrt{w}}\right)^jI_{g+j-1}\left(2\sqrt{w}\right).
\IEEEyesnumber
\end{IEEEeqnarray*}

We find from (8) that the regularized hypergeometric function $\tilde{\Phi} _3(-k,g,\varsigma x,wx)=\frac{1}{\Gamma (g)} \Phi _3(-k,g,\varsigma x,wx)$  with $k\in \mathbb{Z}^+$ can be expressed via the sum of  Bessel I functions as
\begin{IEEEeqnarray*}{c}
\tilde{\Phi} _3(-k;g;\varsigma x, wx)=(wx)^{\frac{1-g}{2}}
\sum _{j=0}^k\frac{(-k)_j}{j!}
\left(\frac{\varsigma}{\sqrt{w} }\right)^j
x^{\frac{j}{2}}I_{g+j-1}\left(2\sqrt{wx}\right).
\IEEEyesnumber
\end{IEEEeqnarray*}
It is interesting to note that the function $\tilde{\Phi} _3$ on the left-hand side (LHS) of (9) is not defined for  $-g\in {\mathbb{Z}^*}$, while the expression on the right-hand-side (RHS) is defined, and it specifies $\underset{g \to m, -m\in \mathbb{Z}^*}{\lim}\tilde{ \Phi} _3(-k ;g ;\varsigma x, wx)$.
 
\subsection{Marcum Q function}
The generalized Marcum Q function is defined as [4]-[8], [16]-[17]
\begin{equation}
Q_M (\alpha ,\beta)=\frac{1}{\alpha ^{M-1}}\int _\beta^\infty x^M\mathrm{exp}\left(-\frac{\alpha ^2+x^2}{2}\right) \mathrm{I}_{M-1}(\alpha x)dx.
\end{equation}
This function has many interesting features, see, for instance, [4]-[8], [16]-[17]. In this work, we employ a differentiation formula w.r.t. $\beta$  given in \cite [eq. (20)] {Bry1} as
\begin{IEEEeqnarray*}{c}
\frac{\partial ^{n}Q_M(\alpha,\beta)}{\partial \beta^{n}}=\frac{n!\alpha^{1-M}}{2^n\beta^{n-M+1}}\mathrm{exp}\left(-\frac{\alpha ^2+\beta ^2}{2}\right)\sum_{k=0}^{\lfloor{n/2}\rfloor}\frac{(-2\alpha ^2)^{-k}}{k!(n-2k)!}\sum_{k=\lfloor n/2\rfloor}^{n}\frac{(-2\beta ^2)^{k}}{(n-k)!(2k-n)!}\\
\times \sum_{j=0}^{k-1}{{k-1}\choose{j}}\left(-\frac{\alpha}{\beta}\right)^jI_{M-j-1}\left(\alpha,\beta\right)
\IEEEyesnumber
\end{IEEEeqnarray*}
where $r \choose q$ is a binomial coefficient.

We also use a recurrence formula for the Marcum Q function given in \cite[eq. (4)]{Bry1} as 
\begin{IEEEeqnarray*}{c}
Q_{M+n}(\alpha ,\beta)=
Q_{M}(\alpha ,\beta)+\left(\frac{\beta}{\alpha}\right)^M\mathrm{exp}\left(-\frac{\alpha ^2+\beta ^2}{2}\right)\sum_{k=0}^{n-1}\left(\frac{\beta}{\alpha}\right)^kI_{M+k}\left(\alpha \beta\right), \hspace{3mm}n\in \mathbb{Z}^+,
\IEEEyesnumber
\end{IEEEeqnarray*}
as well as a formula \cite [eq. (11)] {Dri} that gives a relation between two  Marcum Q functions of positive and negative integer orders as
\begin{equation}
Q_M(\alpha ,\beta)+Q_{1-M}(\beta ,\alpha)=1. 
\end{equation}

 A relation formula between the Marcum Q of the integer order $M$ and confluent hypergeometric $\Phi _3$ functions was recently derived in \cite [eq. (9)] {Mor} as
\begin{IEEEeqnarray*}{c}
Q_M(\alpha ,\beta)=\left(\frac{\alpha ^2}{2}\right)^{1-M}\frac{\mathrm{exp}\left(-\frac{\alpha ^2+\beta ^2}{2}\right)}{\Gamma(2-M)}
\Phi _3\left(1;2-M;\frac{\alpha ^2}{2};\frac{\alpha ^2\beta ^2}{4}\right), \hspace{2mm}M<2.
\IEEEyesnumber
\end{IEEEeqnarray*}

The arguments of the Marcum Q function $\alpha$ and $\beta$ may be real or purely imaginary \cite{Bla}, and in the latter case, the Marcum Q function is referred to as the modified Q function \cite {Bla}.  A series expansion  given in \cite [eq. (2.6)] {Andra} in terms of Laguerre polynomials $L_k^M (.)$ as
\begin{equation}
Q_M (\alpha ,\beta)=1-\mathrm{exp}\left(-\frac{a^2}{2}\right)\sum _{k=0}^{\infty}(-1)^n\frac{L_k^{(M-1)}\left(\frac{a^2}{2}\right)}{\Gamma(M+k+1)}\left(\frac{b^2}{2}\right)^{k+M}
\end{equation}
proves that the modified Marcum Q function is real-valued.

\section{Evaluation of Laplace Transform}
We start with the presentation  of two lemmas required for the evaluation of (1).

\emph{Lemma 1}: 
The confluent hypergeometric function $\Psi _2(a;d,a;w,z)$ is expressed in terms of the confluent hypergeometric function $\Phi _3 $ as
\begin{IEEEeqnarray*}{c}
\Psi _2(a;d ;a ;w;z)
=\mathrm{exp}\left(z+w\right)\Phi _3 (d -a;d;-w,w z), \hspace{5mm} -a ,-d \notin {\mathbb{Z}^*}.
\IEEEyesnumber
\end{IEEEeqnarray*}
\emph{Proof}: See Appendix A.

The next lemma extends previously reported results on the relation between the Marcum Q function and hypergeometric function $\Phi _3 (b, g, t,v)$ (valid only for $g\in\mathbb{Z}^+$) [16]-[17] to the case of $g\in (0, \infty)$.

\emph{Lemma 2}:  The  Marcum Q function $Q_M(\alpha ,\beta )$ with $M>-1$ can be expressed via the hypergeometric function $\Phi _3 (.)$  as 
\begin{IEEEeqnarray*}{c}
1-Q_M(\alpha ,\beta )=\left(\frac{\beta ^2}{2}\right)^M\mathrm{exp}\left(-\frac{\alpha ^2+\beta ^2}{2}\right)\frac{\Phi _3\left(1;M+1;\frac{\beta ^2}{2},\frac{\alpha ^2\beta ^2}{4}\right)}{\Gamma (M+1)},
\IEEEyesnumber
\end{IEEEeqnarray*}
and the
hypergeometric function $\Phi _3(b;g;t,v)$ with $b\in \mathbb{Z}^+$ and $g>0$
can be expressed via the Marcum Q function as 
\begin{IEEEeqnarray*}{c}
\Phi _3(b;g;t,v)=\Gamma (g) \mathrm{exp}\left(v/t+t\right)\sum _{j=0}^{2(b-1)} \delta _j(b,g,t,v)t^{1-g+j}\left[1-Q_{g-1-j}\left(\sqrt {2v/t},\sqrt {2t}\right)\right].
\IEEEyesnumber
\end{IEEEeqnarray*}
The parameter $\delta _j(b,c,w,z)$ in (18) is defined as \cite [eq.(39)] {Mor}
\begin{IEEEeqnarray*}{c}
\delta _j(b,g,w,z)=\frac{(-1)^{b-1}z^{b-1-j}}{w^{b-1}\Gamma (b)}\sum _{k=0}^{\lfloor j/2\rfloor}\frac{(-1)^k(b-j+k)_{j-k}
(g-j-1+k)_{j-2k}}{(j-2k)!k!}z^{k}
\IEEEyesnumber
\end{IEEEeqnarray*}
where $\lfloor u\rfloor$ means the integer part of $u$.

\emph{Proof}: See Appendix B.

It is seen that for integer values of $M$, (17) with the help of (13) reduces to (14) derived in \cite{Mor}.

\emph{Corollary 1}: The regularized hypergeometric function $\tilde{\Psi}_2(a;a+n ;a ;w;z)=\frac{1}{\Gamma (a+n)} \Psi _2(a;a+n ;a ;w;z)$ with $ n\in \mathbb{Z}^+$ and $-a, -(a+n)\notin \mathbb{Z} ^*$ can be expressed in terms of Marcum Q function $Q_M(\alpha ,\beta)$ as
\begin{IEEEeqnarray*}{c}
\tilde{\Psi }_2(a ;a+n ;a ;w;z)
=\sum _{j=0}^{2(n-1)} (-w)^{1-a-n+j}\\
\times \delta _j(n,a+n,-w,wz) \left[1-Q_{a+n-j-1} \left(i\sqrt{2z}, i\sqrt{2w}\right)\right]
\IEEEyesnumber
\end{IEEEeqnarray*}
where $i=\sqrt{-1}$. 

\emph{Proof:} Eq. (20) follows directly from (16) and (18)-(19). 

It is seen that the function on the LHS of (20) is not defined for $-a\in \mathbb{Z}^*$ and for $-(a+n)\in \mathbb{Z}^*$, while the expression on the RHS is defined.  The latter formula specifies $\underset{a \to m, -m\in \mathbb{Z}^*}{\lim} \tilde{\Psi} _2(a; a+n ;a ;w;z)$ and $\underset{(a+n) \to m, -m\in \mathbb{Z}^*}{\lim} \tilde{\Psi} _2(a; a+n ;a ;w;z)$, see also the proof in Appendix A.

\emph{Corollary 2}: The regularized hypergeometric function $\tilde{\Psi}_2( a+n; a;a+n ;w;z)=\frac{1}{\Gamma (a)} \Psi _2(a+n;a;a+n ;w;z)$ with $ n\in \mathbb{Z}^+$ and $-a, -(a+n)\notin \mathbb{Z} ^*$ can be expressed via the sum of Bessel I functions  as
\begin{IEEEeqnarray*}{c}
\tilde{\Psi }_2(a+n ;a ;a+n ;w;z)
=\mathrm{exp}\left(w+z\right)(wz)^{\frac{1-a}{2}}\\
\times \sum _{j=0}^n (-1)^j \frac{(-n)_j}{j!}
\left(\frac{w}{z}\right)^{\frac{j}{2}}I_{a+j-1} \left(2\sqrt{wz}\right).
\IEEEyesnumber
\end{IEEEeqnarray*}

\emph{Proof:} Eq. (21) is obtained immediately from (16) and (9). 

Similarly to (20), the expression on the RHS of (21) defines $\underset{a \to m, -m\in \mathbb{Z}^*}{\lim} \tilde{\Psi} _2(a+n; a ;a+n ;w;z)$ and $\underset{(a+n) \to m, -m\in \mathbb{Z}^*}{\lim} \tilde{\Psi} _2(a+n; a ;a+n ;w;z)$.

Then we consider (1). The solution of (1) depends on the relation between the parameters $\mu _1$ and $\mu _2$, and
it can be obtained via the following proposition.

\emph{Proposition 1}:

1. For the arbitrary real $\mu _2$ and $\mu _2=\mu _1$, (1) can be evaluated as
\begin{IEEEeqnarray*}{c}
\label{dint3}
\mathrm{In}(\alpha ,\beta ,c,p,\mu _1,\mu _1)
=\frac{1}{p}\left(\frac{c}{2p}\right)^{\mu_1 -1}
\mathrm{exp}
\left(\frac{ c^2}{4p}
\right)
Q_{\mu _1}\left(\frac{\alpha c}{\sqrt{2p\cdot\tilde{p}}},\beta \sqrt {\frac{2p}{\tilde{p}}}\right)
\IEEEyesnumber
\end{IEEEeqnarray*}
where $\tilde{p}=2p+\alpha ^2$.

2.  For the real $\mu _2>-1$ and $\mu _1=\mu _2+n$ with $n\in \mathbb{Z}^+$, (1) can be evaluated as
\begin{IEEEeqnarray*}{c}
\mathrm{In}(\alpha ,\beta ,c,p,\mu _2+n,\mu _2)=\mathrm{In}(\alpha ,\beta ,c,p,\mu _2,\mu _2)
+\frac{2}{c}\left(-\frac{\alpha ^2}{c}\right)^{-\mu _2 }\mathrm{exp}\left(-\frac{\beta ^2}{2}\right)\sum_{k=0}^{n-1}\sum _{j=0}^{2k}\left(\frac{\beta ^2}{2}\right)^j\\
\times \left(\frac{-\alpha ^2}{\tilde{p}}\right)^{-k+j}\delta _j\left(k+1,\mu _2+k+1,-\frac{\alpha ^2\beta ^2}{2\tilde{p}},\frac{\alpha ^2\beta ^2c^2}{4\tilde{p}^2}\right) \left[1-Q_{\mu _2+k-j}\left(i\frac{c}{\sqrt{\tilde{p}}},i\frac{\alpha \beta}{\sqrt{\tilde{p}}}\right)\right].
\IEEEyesnumber
\end{IEEEeqnarray*}


3. For  $\mu _2\in (0,\infty)$, and $\mu _1=\mu _2+n$, where $n\in \mathbb{Z}^+$,
 (1) can be also evaluated via the sum of finite-range integrals as
\begin{IEEEeqnarray*}{c}
\mathrm{In}(\alpha ,\beta ,c,p,\mu _2+n,\mu _2)=\mathrm{In}(\alpha ,\beta ,c,p,\mu _2,\mu _2)\\
+\frac{2}{\tilde{p}}\left(\frac{\beta ^2c}{2\tilde{p}}\right)^{\mu2 -1}
\frac{\mathrm{exp}\left(-\frac{\beta ^2}{2}+\frac{c^2}{2\tilde{p}}\right)}{\Gamma(\mu _2) } \sum_{j=1}^{n}\frac{(\beta ^2/2)^{j}}{  (j-1)!}\\
\times \int _0^1t^{\mu _2-1}(1-t)^{j-1}\mathrm{exp}\left(\frac{\alpha ^2\beta ^2}{2\tilde{p}}t\right){_0F_1}\left[;\mu _2;\frac{\alpha ^2\beta ^2c^2}{4\tilde{p}^2}t\right]dt
\IEEEyesnumber
\end{IEEEeqnarray*}
where $_0F_1(;q,z)=\Gamma (q)z^{(1-q)/2}I_{q-1}\left(2\sqrt{z}\right)$ is a hypergeometric function \cite [vol. 3, eq. (7.13.1.1)] {Prud}.


4.   A solution to (1) for  $\mu _1=\mu _2-n$, where $n\in \mathbb{Z}^+$,
can be expressed as
\begin{IEEEeqnarray*}{c}
\mathrm{In}(\alpha ,\beta ,c,p,\mu _1,\mu _1+n)
=\mathrm{In}(\alpha ,\beta ,c,p,\mu _1+n,\mu _1+n)\\
-\frac{2}{c}
\left(\frac{c}{\tilde{p}}\right)^{n}\left(\frac{\beta}{\alpha}\right)^{\mu _1}\mathrm{exp}\left(\frac{c^2/2-p\beta ^2}{\tilde{p}}\right)\sum_{k=0}^{n-1}\sum_{j=0}^{n-k-1}\left(\frac{\beta }{c}\right)^{j+k}\\
\times\left(\frac{\tilde{p}}{\alpha}\right)^k\left(-\alpha\right)^j \frac{(-n+k+1)_j}{j!}I_{\mu _1+k+j}\left(c\frac{\alpha \beta }{\tilde{p}}\right).
\IEEEyesnumber
\end{IEEEeqnarray*}

\emph{Proof}: See Appendix C.


\emph{Corollary 3}: Proposition 1 gives also solutions to the integral having the form of \\$\mathcal{L}\left\{Q_{\mu _1}(\alpha ,\beta \sqrt{t})t^{\frac{\mu _2-1}{2}}I_{\mu_2 -1}(c\sqrt{t}),\left\{t,p\right\}\right\}$, where $\mu _1$ and $\mu _2$  are integers. It can be expressed as
\begin{IEEEeqnarray*}{c}
\mathcal{L}\left\{Q_{\mu _1}(\alpha ,\beta \sqrt{t})t^{\frac{\mu _2-1}{2}}I_{\mu_2 -1}(c\sqrt{t}),\left\{t,p\right\}\right\}\\
=\left(\frac{c}{2}\right)^{\mu _2-1}p^{-\mu _2}\mathrm{exp}\left(\frac{c^2}{4p}\right)-
\mathrm{In}(\beta ,\alpha ,c,p,1-\mu _1,\mu _2).
\IEEEyesnumber
\end{IEEEeqnarray*}

\emph{Proof}: Eq. (26) follows immediately from  (13) \cite [eq. (11)] {Dri},
and a  Laplace transform formula \cite [vol. 4, eq. (3.15.2.8)] {Prud}
\begin{IEEEeqnarray*}{c}
\mathcal{L}\left\{t^{\frac{\mu -1}{2}}I_{\mu -1}(c\sqrt{t}),\left\{t,p\right\}\right\}=\left(\frac{c}{2}\right)^{\mu -1}p^{-\mu }\mathrm{exp}\left(\frac{c^2}{4p}\right).
\IEEEyesnumber
\end{IEEEeqnarray*} 

We see that  Proposition 1 gives alternative solutions for the cases of $\mu _1=\mu _2+n$ with $n\in\mathbb{Z}^+$ either via the sum of modified Marcum Q functions (24) or via the sum of finite-range integrals (25). The finite-range integral in (25) is proper, and it can be easily evaluated via the standard software. 
\section{ Applications}
In this section, we present  examples of application of the derived results to the performance evaluation of communication systems operating over $\kappa$-$\mu$ fading.
\subsection { Outage probability  analysis in $\kappa$-$\mu$/$\kappa$-$\mu$ interference-limited scenarios}
The problem of analyzing communication systems under co-channel interference (CCI) arises in many practical applications such as cellular, and  ad hoc networks.  This problem is very important in cognitive radio systems.  The analysis of outage probability (OP) under CCI was presented  for different signal of interest (SoI) and CCI fading models, see, for example, [23]-[25] and [13]-[15]. 

The OP in interference-limited scenarios is merely the CDF of a random variable (RV) represented by the ratio of SoI and CCI powers, which is defined on the basis of a well-known rule \cite{Rom}, \cite{Pap} as 
\begin{IEEEeqnarray*}{c}
P_{\mathrm{out}}(z)=\mathrm{Pr}\left \{\frac {\mathrm{Power}_{\mathrm{SoI}}}{\mathrm{Power}_{\mathrm{CCI}}}<z\right \}
=\int _{0}^{\infty}F_{\mathrm{SoI}}\left (z y \right)f_{\mathrm{CCI}}(y)dy
\IEEEyesnumber
\end{IEEEeqnarray*}
where $f_x(z)$ and  $F_x(z)$ are the respective PDF and  CDF of the RV $x$.

Let  the SoI statistical parameters be $\kappa _{\mathrm{S}}$, $\mu _{\mathrm{S}}$, and $\Omega _{\mathrm{S}}$, while the CCI parameters be $\kappa _{\mathrm{I}}$, $\mu _{\mathrm{I}}$, and $\Omega _{\mathrm{I}}$. Then plugging (2)-(3) into (28) and using (13), we find that for integer values of $\mu _{\mathrm{S}}$ and $\mu _{\mathrm{I}}$, the OP can be expressed as
\begin{IEEEeqnarray*}{c}
P_{\mathrm{out}}(z)=
\frac{\mu _{\mathrm{I}} (1+\kappa _{\mathrm{I}})^{\frac{\mu _{\mathrm{I}}+1}{2}}} 
{\kappa _{\mathrm{I}}^{\frac{\mu _{\mathrm{I}}-1}{2}}\mathrm{exp}\left(\mu _{\mathrm{I}} \kappa _{\mathrm{I}}\right)\Omega _{\mathrm{I}}^{\frac{\mu _{\mathrm{I}}+1}{2}} }
\mathrm{In}\left(\beta ' ,\alpha ' ,c',p',1-\mu _{\mathrm{S}},\mu _{\mathrm{I}}\right)
\IEEEyesnumber
\end{IEEEeqnarray*}
where $c'=2 \mu _{\mathrm{I}}\sqrt{\frac{\kappa _{\mathrm{I}} (1+\kappa _{\mathrm{I}})}{ \Omega _{\mathrm{I}}}}$, $p'=\frac{ \mu _{\mathrm{I}} (1+\kappa _{\mathrm{I}}) }{\Omega _{\mathrm{I}}}$, $\alpha '=\sqrt{2\kappa _{\mathrm{S}}\mu _{\mathrm{S}}}$, and $\beta'=\sqrt{\frac{2(1+\kappa _{\mathrm{S}})\mu _{\mathrm{S}}z}{\Omega _{\mathrm{S}}}}$.

In Fig. 1, we present numerical estimates of the OP depicted versus the signal-to-interference ratio (SIR) defined as $\Omega _{\mathrm{S}}/\Omega _{\mathrm{I}}$. The estimates in Fig. 1 are given for single-input single-output scenarios with $\kappa _{\mathrm{S}}=0.5$ and $\kappa _{\mathrm{S}}=2.5$; $\mu _{\mathrm{S}}=\mu _{\mathrm{I}}=2$ and $\kappa _{\mathrm{I}}=0.5$.
 \subsection {Probability of energy detection of an unknown signal in  $\kappa$-$\mu$ fading channels}
The detection of unknown signals is an important issue in various applications such as cognitive radio systems \cite{Sim2}, [11]-[12]. The detection procedure can be implemented in different ways. Measuring the energy of the received waveform over an observation time window and then comparing it with a threshold, is the simplest detection method referred to as the energy detection \cite{Sim2}. Evaluation of the detection probability, $P_{\mathrm{d}}$, over $\kappa$-$\mu$ fading was solved so far by applying series expansions of Marcum Q and Bessel I functions [11]-[12]. These procedures resulted in formulas given in terms of infinite series \cite{Annam} that can be expressed in terms of multivariate hypergeometric functions \cite{Sof}. The hypergeometric functions required for the evaluation are not nowadays implemented in a standard software. Meanwhile, closed-form expressions can be obtained for some combinations of parameters  based on the results of this work.

   Respective expressions for the probability of detection and false alarm ($P_{\mathrm{f}}$) in an AWGN channel are given in  \cite {Sim2}, [11]-[12] as
\begin{IEEEeqnarray*}{c}
P_{\mathrm{d}}=Q_u\left(\sqrt{2\gamma},\sqrt{2\lambda}\right),
\IEEEyesnumber
\end{IEEEeqnarray*}
and 
\begin{IEEEeqnarray*}{c}
P_{\mathrm{f}}=\frac{\Gamma\left(u,\lambda /2\right)}{\Gamma\left(u\right)}
\IEEEyesnumber
\end{IEEEeqnarray*}
where $u$ is the product of the observation time and signal bandwidth, $\Gamma\left(a,x\right)=\int _x^\infty t^{a-1}\mathrm{exp}\left(-t\right)dt$ is the upper incomplete gamma function, $\gamma$ denotes the signal-to-noise ratio (SNR), and $\lambda$ is the detector threshold.
Thus, by averaging (30) over the fading distribution (2), we obtain that the probability of energy detection over $\kappa$-$\mu$ fading,$P_{\mathrm{d}_{\kappa-\mu }}$,  can be expressed as
\begin{IEEEeqnarray*}{c}
P_{\mathrm{d}_{\kappa-\mu }}=
\frac{\mu (1+\kappa )^{\frac{\mu +1}{2}}} 
{\kappa ^{\frac{\mu -1}{2}}\mathrm{exp}\left(\mu \kappa\right)\Omega ^{\frac{\mu +1}{2}} }\mathrm{In}\left(\sqrt{2},\sqrt{2\lambda}, \hat{c},\hat{p},u,\mu \right)
\IEEEyesnumber
\end{IEEEeqnarray*}
where $\hat{c}=2 \mu \sqrt{\frac{\kappa (1+\kappa )}{ \Omega }}$ and $\hat{p}=\frac{ \mu  (1+\kappa )  }{\Omega }$.

In Fig. 2, we present average estimates of the detection probability over $\kappa$-$\mu$ fading for $u=2.5$, $\mu =0.5$, and for two values of $\kappa$: $\kappa =0.5$ and $\kappa =4.2$.  
\begin{figure}[!t] 
\centering
\includegraphics[width=4.2in,height=4in]{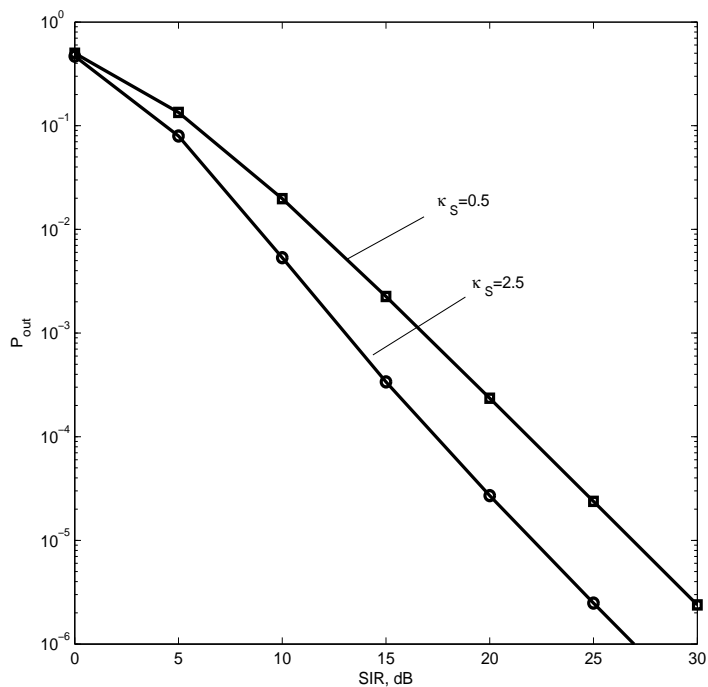}
\vspace{-5mm}
\caption{Outage probability versus the signal-to-interference ratio in $\kappa$-$\mu$/$\kappa$-$\mu$ 
interference-limited scenarios. Single points report simulation results.}
\label{fig:1}
\end{figure}  
\begin{figure}[!t] 
\centering
\includegraphics[width=4.2in,height=4in]{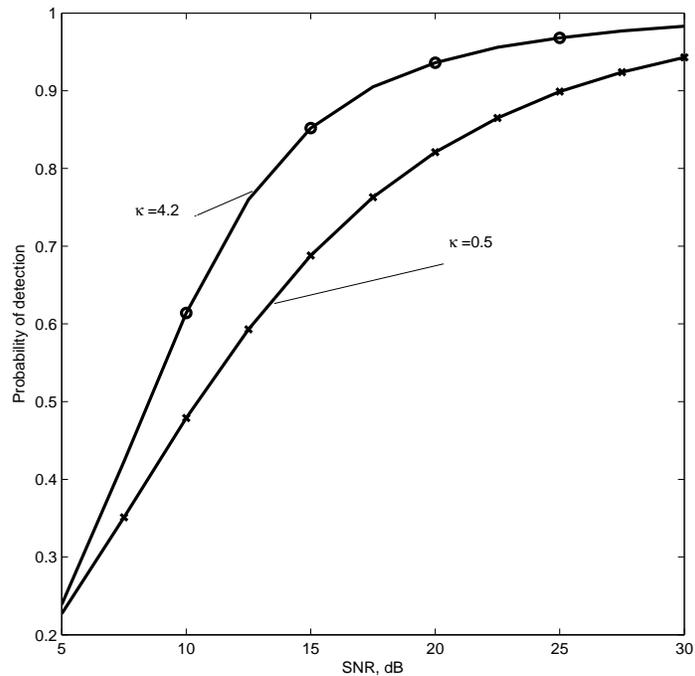}
\vspace{-5mm}
\caption{Probability of detection  in $\kappa$-$\mu$ fading, $P_{\mathrm{f}}=0.1$. Single points report simulation results.}
\label{fig:1}
\end{figure}  
\section{Conclusion}
In this paper, we solve the infinite integral (1) that is required for assessing a few performance metrics of communication systems operating over $\kappa$-$\mu$ fading channels. On the basis of derived results, the  probability of  detection of unknown signals and outage probability in $\kappa$-$\mu$/$\kappa$-$\mu$ interference-limited scenarios can be evaluated analytically. 

 Since the $\kappa$-$\mu$ fading distribution is a general statistical model comprising many fading distributions, the results of this work can be applied to a large variety of fading scenarios. The analytical results derived in this paper are restricted by scenarios with $\mu _2=\mu _1+n$, where $n$ is an integer, see Proposition 1. Based on the fact that the parameter $\mu$ of the $\kappa$-$\mu$ distribution is inversely proportional to the amount of fading under a fixed $\kappa$ \cite [eq. (8)] {Yac}, our results can be used as bounds for real estimates in scenarios where $n \notin \mathbb{Z}$. 

All derived formulas are given in terms of  functions available via the  standard software such as \emph{Mathematica}. The theoretical contribution of this work includes also the recognition of relations between the hypergeometric functions of two variables $\Psi _2$ and $\Phi _3$, see (16), as well as an extension of previously derived results on the connection between the generalized Marcum Q function and confluent hypergeometric function $\Phi _3$, see (17)-(18). These facts resulted in establishing a connection between the hypergeometric function $\Psi _2$ and modified  Marcum Q function, see (20), as well as between $\Psi _2$ and  Bessel I functions, see (21).

\appendices
\section{Proof of Lemma 1}
We use a series expansion of $\Psi _2(a ;d ,a ;w,z)$ given in \cite [vol. 3, eq. (6.6.2.2)] {Prud} as
\begin{IEEEeqnarray*}{c}
\Psi _2(a ;d, a ;w,z)=\sum _{k=0}^{\infty}\frac{(wz)^k}{(b )_kk!}{_1F_1\left(a +k;a +k;z\right)}{_1F_1\left(a +k;d +k;w\right)}\\
=\mathrm{exp}(z)\sum _{k=0}^{\infty}\frac{(wz)^k}{(d )_kk!}{_1F_1\left(a +k;d +k;w\right)}
\IEEEyesnumber
\end{IEEEeqnarray*}
where \cite [vol. 3, eq. (7.2.2.1)]{Prud} and \cite [vol. 3, eq. (7.2.2.8)]{Prud}
\begin{IEEEeqnarray*}{c}
{_1F_1(s;q;w)}=\sum _{l=0}^{\infty}\frac{(s)_l}{(q)_ll!}w^l =\mathrm{exp}(w){_1F_1}\left(q-s;q;-w\right),\hspace{5mm}
-q\notin \mathbb{Z}^{*},
\IEEEyesnumber
\end{IEEEeqnarray*}
 is a confluent hypergeometric function. By plugging (35) into (34) we obtain the product of $\mathrm{exp}(w+z)$ and the hypergeometric series (5). \hspace{80mm}                                                                                               \IEEEQED

\section{Proof of Lemma 2}

We evaluate the Laplace transforms of functions on the LHS and RHS of (17) w.r.t. $\frac{\beta ^2}{2}$.
 Using \cite[vol. 4, eq. 1.1.3.1]{Prud} and \cite[vol. 4, eq. 3.43.8]{Prud} (valid for $M>-1$) we find that
\begin{IEEEeqnarray*}{c}
\mathcal{L}\left\{\mathrm{exp}\left(-\frac{\alpha ^2+\beta ^2}{2}\right)\left(\frac{\beta ^2}{2}\right)^M\frac{\Phi _3\left(1,M+1,\frac{\beta ^2}{2},\frac{\alpha ^2\beta ^2}{4}\right)}{\Gamma (M+1)},\{\beta ^2/2,p\}\right\}\\
=\frac{\mathrm{exp}\left(-\frac{\alpha ^2}{2}\right)}{p(p+1)^M}\mathrm{exp}\left[\frac{\alpha ^2}{2(p+1)}\right].
\IEEEyesnumber
\end{IEEEeqnarray*}
Evaluating $\mathcal{L}\left\{1-Q_M(\alpha,\beta,\{\beta ^2/2,p\}\right\}$ we obtain that
\begin{IEEEeqnarray*}{c}
\mathcal{L}\left\{1-Q_M(\alpha,\beta),\{\beta ^2/2,p\}\right\}=\frac{\mathrm{exp}\left(-\frac{\alpha ^2}{2}\right)2^{\frac{M-1}{2}}}{\alpha ^{M-1}}\\
\times \mathcal{L}\left\{\int _0^{\frac{\beta ^2}{2}} \mathrm{exp}\left(-t\right)t^{\frac{M-1}{2}}I_{M-1}(\alpha \sqrt{2t})dt,\{\beta ^2/2,p\}\right\}\\
=\frac{\mathrm{exp}\left(-\frac{\alpha ^2}{2}\right)}{p(p+1)^M}\mathrm{exp}\left[\frac{\alpha ^2}{2(p+1)}\right]
\IEEEyesnumber
\end{IEEEeqnarray*}
where we used (10) and Laplace transform properties given by \cite[vol. 4, eq. 1.1.3.1]{Prud}, \cite[vol. 4, eq. 1.1.5.2]{Prud}, and \cite[vol. 4, eq. 3.15.2.8]{Prud}.

We see that the expressions on the RHS of (35) and (36) are equal, and thus (17) holds true. 
 Eq. (18) for $b=1$ is obtained directly from (17). Then applying a recurrence method proposed in \cite[eq. (36)-(37)]{Mor} and increasing $b$, we obtain (18) with the parameter $\delta _j(b, g, w,z)$ equal to that derived in \cite[eq. (39)]{Mor}. In such a way we express $\tilde{\Phi} _3$ directly in terms of Marcum Q function.  \hspace{145mm} \IEEEQED
\section{Proof of Proposition 1}

Using the differentiation formula (11) \cite [eq. (20)] {Bry1} for $n=1$ we obtain that
\begin{IEEEeqnarray*}{c}
\label{dint1}
\frac{\partial \mathrm{In_1}(\alpha ,\beta ,c,p,\mu _1,\mu _2)}{\partial \beta}=-\frac{2\beta ^{\mu _1}}{\alpha ^{\mu _1-1}}\mathrm{exp}\left(-\frac{\beta ^2}{2}\right)
\int_0^{\infty}\mathrm{exp}\left[-\left(p+\frac{\alpha ^2}{2}\right)t^2\right]t^{\mu _2-\mu _1+1}\\
\times I_{\mu _1-1}(\alpha \beta t)I_{\mu_2 -1}
(ct)dt.
\IEEEyesnumber
\end{IEEEeqnarray*}

If $\mu _1=\mu _2$, (\ref {dint1}) can be evaluated directly on the basis of \cite[vol. (4), eq. (3.15.17.1)]{Prud}, resulting in  
\begin{equation}
\label{dint2}
\frac{\partial \mathrm{In_1}(\alpha ,\beta ,c,p,\mu _1,\mu _1)}{\partial \beta}=-\frac{2\beta ^{\mu _1}}{\alpha ^{\mu _1-1}\left(2p+\alpha ^2
\right)}\mathrm{exp}
\left[\frac{-2p\cdot \beta ^2+c^2}{2(2p+\alpha ^2)}\right]
I_{\mu _1-1}\left(\frac{\alpha \beta c}{2p+\alpha ^2}\right),
\end{equation}
and (22)  follows directly from (38) and (10).  

 From (12) \cite [eq. (4)] {Bry1},
we obtain  a recurrence relation for (1) with $\mu _1=\mu _2+n$ as
\begin{IEEEeqnarray*}{c}
\mathrm{In}(\alpha ,\beta ,c,p,\mu _2+n,\mu _2)=\mathrm{In}(\alpha ,\beta ,c,p,\mu _2,\mu _2)+\mathrm{exp}\left(-\frac{\beta ^2}{2}\right)
\left(\frac{\beta }{\alpha }\right)^{\mu _2}\\
\times \sum_{k=0}^{n-1}\left(\frac{\beta }{\alpha }\right)^{k}\int_0^{\infty}\mathrm{exp}\left[-(p+\alpha ^2/2)t\right]t^{\frac{-k-1}{2}}I_{\mu _2+k}(\alpha \beta \sqrt{t})I_{\mu_2 -1}(c\sqrt{t})dt.
\IEEEyesnumber
\end{IEEEeqnarray*}
Applying a Laplace transform formula \cite [vol. 4, eq. (3.15.17.13)] {Prud}, we find that
\begin{IEEEeqnarray*}{c}
\mathrm{In}(\alpha ,\beta ,c,p,\mu _2+n,\mu _2)=\mathrm{In}(\alpha ,\beta ,c,p,\mu _2,\mu _2)+2c^{\mu _2-1}\mathrm{exp}\left(-\frac{\beta ^2}{2}\right)
\left(\frac{\beta ^2}{4p+2\alpha ^2 }\right)^{\mu _2}\\
\times \sum_{k=0}^{n-1}\frac{\left(\beta ^2/2\right)^k}{\Gamma (\mu _2+k+1) }\Psi _2\left[\mu _2;\mu _2+k+1, \mu _2; \frac{\alpha ^2\beta ^2}{4p+2\alpha ^2},\frac{c ^2}{4p+2\alpha ^2}\right].
\IEEEyesnumber
\end{IEEEeqnarray*}
Then using Corollary 1, we obtain (23).

Eq. (24) follows from (40), (33), and an integral representation of ${_1F_1}$ in (33), which is given in \cite [vol. 3, eq. (7.2.2.6)] {Prud} as 
\begin{IEEEeqnarray*}{c}
{_1F_1(a;b;z)}=\frac{\Gamma(b)}{\Gamma (a)\Gamma (b-a)}\int _0^1t^{a-1}(1-t)^{b-a-1}\mathrm{exp}(zt)dt,\hspace{10mm} b>a>0.
\IEEEyesnumber
\end{IEEEeqnarray*}
The absolute convergence of $\Psi _2$ justifies changing the order of summation and integration in (33) after using (41). This procedure results in (24).

 For $\mu _1=\mu _2-n$, we find from (12) that
\begin{IEEEeqnarray*}{c}
\mathrm{In}(\alpha ,\beta ,c,p,\mu _1,\mu _1+n)=\mathrm{In}(\alpha ,\beta ,c,p,\mu _1+n,\mu _1+n)
-\mathrm{exp}\left(-\frac{\beta ^2}{2}\right)\sum _{k=0}^{n-1}\left(\frac{\beta}{\alpha}\right)^{\mu _1+k}\\
\int _0^{\infty}\mathrm{exp}\left[-(p+\alpha ^2/2)t\right]t^{\frac{n-k-1}{2}}I_{\mu _1+k}\left(\alpha \beta \sqrt{t}\right)I_{\mu _2-1}\left(c \sqrt{t}\right)dt\\
=\mathrm{In}(\alpha ,\beta ,c,p,\mu _1+n,\mu _1+n)-\frac{2}{c}\left(\frac{c}{\tilde{p}}\right)^{\mu _1+n}\mathrm{exp}\left(-\frac{\beta ^2}{2}\right)\\
\times \sum _{k=0}^{n-1}\left(\frac{\beta ^2}{2}\right)^k \tilde{\Psi }_2\left(\mu _1+n;\mu _1+k+1,\mu _1+n;\frac{\alpha ^2\beta ^2}{4p+2\alpha ^2}, \frac{c ^2}{4p+2\alpha ^2}
\right)
\IEEEyesnumber
\end{IEEEeqnarray*}
where we again used  \cite [vol. 4, eq. (3.15.17.13)] {Prud}. Then applying Corollary 2, we obtain (25).\IEEEQED

\bibliographystyle{IEEEtran}

\end{document}